# Revisiting a Nice Cycle Lemma and its Consequences

M.O. Albertson (posthumous)
Smith College, Northampton, MA, USA
J.P. Hutchinson (emerita)
Macalester College, St. Paul, MN, USA
hutchinson@macalester.edu
R.B. Richter
Department of Combinatorics and Optimization
Faculty of Mathematics, University of Waterloo
Waterloo, ON, CANADA
brichter@uwaterloo.ca

DRAFT - Apr. 10, 2016

With thanks to Xingxing Yu

**Abstract.**
We correct some errors and omissions primarily in the paper [AH04], discovered by the third author, and also some in a proof of [T93] and of [Y97]. We give a short proof of Thomassen's theorem [T93] that every triangulation of a surface with all noncontractible cycles sufficiently long can be 5-colored; part of the shortness is due to the use of the Four Color Theorem, which is not used in Thomassen's original proof.

**0. Introduction**

In this paper let $G$, a simple graph, always have a 2-cell triangulation embedding $\Pi$ of a nonplanar surface $\Sigma$, and suppose the embedding has *width* $w(G, \Pi)$, the length of the shortest noncontractible cycle in $G$, to be determined when and as needed. Such an embedded graph is called *locally planar* when the width is sufficiently large (for the problem at hand.) We prove Thomassen's theorem that locally planar graphs are 5-colorable; it is also known that locally planar graphs are 5-choosable by a much more complex proof [DKM08]. Our proof is shorter than that of [T93] due to the introduction of the idea of *nice cycles* in embedded graphs and the use of the Four Color Theorem [ApHa76, RSST97].

**1. Planarizing cycles, background and new definitions**

In this paper let the *length* of a path (or cycle) of a graph $G$ be the number of edges. For $P$ a path, we denote its length by $l(P)$ and set $|P| = |V(P)|$.

If $U \subset V(G)$, set $N^0(U) = U$, and let $N^1(U) = N(U)$ denote the set of vertices that are not in $U$ but are adjacent to at least one vertex in $U$. Inductively for $i \geq 1$ $N^i(U)$ denotes the set of vertices that are in neither $N^{i-1}(U)$ nor (when $i \geq 2$) in $N^{i-2}(U)$, but are adjacent to at least one vertex of $N^{i-1}(U)$. Thus $N^i(U)$ consists of those vertices of $G$ whose distance from $U$ is exactly $i$.

Given a traversal of $C$, a 2-sided (or orientation preserving) cycle within $G$ with triangulation embedding $\Pi$, we define $R(C)$ and $L(C)$, respectively, the (not necessarily disjoint) right and left neighbors of $C$. $R(C)$ (resp., $L(C)$) is the set of vertices of $N(C)$ that are adjacent to $C$ along an edge emanating from the right (resp., left) side of $C$. We may also say that the edges joining a vertex of $C$ with one of $R(C)$ (resp., $L(C)$) lie on the right (resp., left) of $C$. All separating cycles of $G$ are 2-sided as are all noncontractible cycles when $\Sigma$ is orientable. For $i \geq 1$ we define $R^i(C)$ to be the vertices of $N^i(C)$ that can be reached by a shortest path starting in $C$ with an edge to a vertex of $R(C)$; we define $L^i(C)$ in an analogous way.

A nonorientable surface may also contain 1-sided (or orientation reversing) noncontractible, necessarily nonseparating, cycles. When $C$ is 1-sided, we define $R$ and $L$ locally along $M$, a subpath of $m \geq 1$ edges of $C$ when $C$ has length $k \geq 2\,m$. On a traversal of $M \subset C$ there are (not necessarily disjoint) right neighbors $R(M)$ and left neighbors $L(M)$, both subsets of $N(C)$. Formally $R(M) = \{y \in N(C) :$ there is a vertex $x \in M$ and an edge $xy$ lying on the right side of $M$ in the traversal$\}$. We may say that locally edge $xy$ lies on the right of $C$. $L(M)$ is defined in the same manner.

(Since in the following we consider only chordless noncontractible cycles, we need not be concerned whether cycle-chords lie to the left or to the right of a noncontractible cycle.)

For $M$ a subpath of chordless cycle $C$ (when $C$ is 1- or 2-sided) and for a vertex $x \in M$, we define $R(x)$ to be its right neighbors with respect to a traversal of $M$, and $\deg_R(x) = |R(x)|$ is called the *right degree of $x$*. We make the analogous definition for the *left degree.*

By *cutting* along a chordless, noncontractible cycle $C$, we mean the following transformation of $G$ and its 2-cell embedding $\Pi$ on a surface $\Sigma$ of Euler genus $g > 0$. If $C$ is 2-sided, we retain two copies of $C$, called $C'$ and $C$. Given an orientation of the original $C$, we connect every vertex of $C'$ to vertices of $R(C)$ by edges lying on the right side of $C$, and we connect every vertex of $C$ to $L(C)$ by edges lying on the left side $C$. When $C = \{x_1, \ldots, x_k\}$ is 1-sided, we replace $C$ by $C'$ a cycle of $2k$ vertices $\{x_1, \ldots, x_k, x'_1, \ldots, x'_k\}$, with each vertex $x_i$ replaced by two copies, $x_i$ and $x'_i$, so that, for $i = 1, \ldots k$, $x_i$ is adjacent to each vertex of $L(x_i)$, with $L$ locally defined, and each $x'_i$ is adjacent to each vertex of $R(x_i)$, locally defined. Since (as assumed) $\Pi$ is a triangulation of $\Sigma$, we think of sewing in a disk (or two) to span $C'$ (and $C$), we add a vertex $x'$ adjacent to every vertex of $C'$ and (when $C$ is 2-sided) add a vertex $x$ adjacent to every vertex of $C$. This graph is denoted $G'$ with a 2-cell triangulation embedding $\Pi'$, embedded on a surface (or the union of two surfaces) $\Sigma'$. When $C$ is 1-sided, then the genus of $\Sigma'$ is $g - 1$. When $C$ is 2-sided and nonseparating, then the genus of $\Sigma'$ is $g - 2$. Otherwise when $C$ is 2-sided and separating, $\Sigma'$ is the union of two surfaces with genus defined to be the maximum of the two genera, which sum to $g$. Similarly $w(G', \Pi')$ is defined to be the minimum of the width of the graph on the two surfaces. (See [T93, MT01].)

Lemma 1 [T93, Y97]. Suppose $G$, $\Pi$ is a triangulation embedding on $\Sigma$ of Euler genus $g > 0$. If $C$ is a shortest noncontractible cycle, then $C$ is chordless. In addition if cutting along $C$ gives triangulation $G'$, $\Pi'$ on $\Sigma'$, then $w(G', \Pi') \geq w(G, \Pi)/2$.

The next distance parameter for 2-sided cycles is alluded to in [T93].

Definition. If $C$ is a 2-sided, chordless, noncontractible cycle in an embedded graph $G$, let $\text{ctdist}_G(C)$ (or $\text{ctdist}(C)$ when the graph context is clear) denote the length of a shortest path (or cycle) beginning with an edge to a vertex on one side of $C$ and ending with an edge on the other side of $C$. If there is no such path, declare $\text{ctdist}(C)$ to be infinite.

Note that if $C$ is a 2-sided, chordless, noncontractible cycle in $G$, $\Pi$, and $P$ is a path that defines $\text{ctdist}_G(C)$, then when $P$ is a cycle, it is necessarily noncontractible. When $P$ begins at $x \in C$ and ends at $y \in C$, $x \neq y$, then both $P \cup C_{x,y}$ and $P \cup C_{y,x}$ are necessarily noncontractible where we define $C_{x,y}$ to be the shorter subpath of $C$ joining $x$ and $y$ and $C_{y,x}$ the longer such subpath; see [T93]. We could define $\text{ctdist}(C)$ for $C$ 1-sided also, but then the latter property need not hold. Instead we focus on a variation of this distance parameter.

Definition. Let $C$ be a chordless, noncontractible cycle in $G$, $\Pi$. Let $\text{ncdist}_G(C)$ (or $\text{ncdist}(C)$ when the graph context is clear) denote the length of a shortest path $P$ that begins at $x \in C$ and ends at $y \in C$ such that when $x \neq y$ $P \cup C_{x,y}$ and $P \cup C_{y,x}$ are both noncontractible or when $x = y$ $P$ is noncontractible. If there is no such path, declare $\text{ncdist}(C)$ to be infinite.

Note that for $C$ 2-sided, $\text{ncdist}(C) \leq \text{ctdist}(C)$. When $G$ is embedded on the torus, $\text{ncdist} = \text{ctdist}$. For other surfaces a path $P$ that gives ncdist may begin and end on the same side of $C$ so that possibly ncdist < ctdist.

Proposition 1. Suppose $C$ is a shortest noncontractible cycle in an embedded graph $G$, $\Pi$ of width $w(G, \Pi)$. Then ncdist$(C) \geq w(G, \Pi)/2$.

Proof. By Lemma 1, $C$ is chordless. Let $P$ be a path (or cycle) defining ncdist$(C)$. If $P$ begins and ends at the same vertex, then the cycle $P$ is noncontractible, and $l(P) \geq w(G, \Pi)$. Otherwise if $P$ begins at $x \in C$ and ends at $y \in C$, $y \neq x$, let $C_{x,y}$ denote the shorter subpath of $C$ from $x$ to $y$. Since $l(C_{x,y}) \leq w(G, \Pi)/2$, and the noncontractible cycle $P \cup C_{x,y}$ has length at least $w(G, \Pi)$, we have $l(P) \geq w(G, \Pi)/2$. QED

When $C$ is a contractible cycle in $G$, $\Pi$, we denote by Int$(C)$ the subgraph contained within its interior together with $C$.

Lemma 2 [M92, T93, Y97]. For $k \geq 1$, let $\Pi$ be a triangulation embedding of $G$ on $\Sigma$ of Euler genus $g \geq 1$ with $w(G, \Pi) \geq 2k+2$. Let $x \in V(G)$. Then $G$ has a unique chordless $\Pi$-contractible cycle $D$ with $x$ in its interior such that every vertex of $D$ has distance $k$ from $x$, and for each $z \in D$ there is a path in Int$(D)$ of length $k$ from $x$ to $z$.

When $X \subset V(G)$, let $G[X]$ denote the induced subgraph on the vertices of $X$. Then the cycle $D$ of Lemma 2, called the *k-th canonical cycle of $x$* in $G$, is a subgraph of $G[N^k(x)]$.

Corollary 1. For $k \geq 1$, let $\Pi$ be a triangulation embedding of $G$ on $\Sigma$ of Euler genus $g \geq 1$ with $w(G, \Pi) \geq 2k+2$. After cutting along a 2-sided, chordless, noncontractible cycle $C$, let the resulting triangulation be $G'$ with triangulation $\Pi'$ on $\Sigma'$ with added vertices $x'$ and $x$. Then if $D'$ is the $k$-th canonical cycle of $x'$ in $G'$ and $D$ the $k$-th canonical cycle of $x$ in $G'$, then Int$(D')$ and Int$(D)$ are disjoint and $D'$ and $D$ are joined by no edge. In particular $R^k(C) \cap L^k(C) = \emptyset$ and no edge joins $R^k(C)$ with $L^k(C)$.

Proposition 2. For $k \geq 1$, let $\Pi$ be a triangulation embedding of $G$ on $\Sigma$ of Euler genus $g \geq 1$ with $w(G, \Pi) \geq 2k+4$. Let $C$ be a 1-sided, chordless, noncontractible cycle $C$, and let $M$ be a subpath of $C$ with at most $m \leq |C|/2$ edges. If $P$ is a path starting from a vertex $x \in M$ on an edge to $R(M)$ and ending at $y \in M$ on an edge from $L(M)$, then $l(P) \geq |C|/2$. In particular $R^{k/2}(M) \cap L^{k/2}(M) = \emptyset$ and no edge joins $R^{k/2}(M)$ with $L^{k/2}(M)$.

Proof. Let $C_{x,y}$ denote the shorter subpath of $C$ between $x$ and $y$ and $C_{y,x}$ the longer. At least one of $P \cup C_{x,y}$ and $P \cup C_{y,x}$ is noncontractible, and $C_{x,y} \subseteq M$. If $P \cup C_{x,y}$ is noncontractible, then as argued in the proof of Prop. 1, $l(P) \geq |C|/2$. Then since $l(P) \geq k+2$, $R^{k/2}(M)$ and $L^{k/2}(M)$ are disjoint and joined by no edge. Otherwise $P \cup C_{x,y}$ is contractible so that $P$ begins and ends with two edges on the same side of $C_{x,y}$, a contradiction. QED

For $G$ with embedding $\Pi$ in $\Sigma$ of Euler genus $g > 0$, a set of pairwise disjoint, noncontractible cycles $C_1, C_2, ..., C_s$ is called a planarizing set if $G - \cup V(C_i)$ is a planar graph. For $i \neq j$, $dist(C_i, C_j)$ is the length of a shortest path joining these two cycles.

We (the second two authors) reprove the following theorem to encorporate part (iii) (and to correct some errors and omissions in earlier proofs.)

Theorem 1 (See also [T93, Y97, MT01, AH04]). For $d \geq 1$, let $\Pi$ be a 2-cell triangulation embedding of a simple graph $G$ on a surface $\Sigma$ of Euler genus $g \geq 1$. Suppose $w(G, \Pi) \geq 8(d+1)(2^{g/2} - 1)$ when $\Sigma$ is orientable and $w(G, \Pi) \geq 8(d+1)(2^g - 1)$ otherwise. Then $G$ contains a planarizing set of cycles $S_1, S_2, ..., S_k$ with $k = g/2$ when $\Sigma$ is orientable and with $g/2 \leq k \leq g$ otherwise, such that
    (i) each $S_i$ is chordless, $1 \leq i \leq k$,
    (ii) dist$(S_i, S_j) \geq d$, $1 \leq i < j \leq k$, and
    (iii) ncdist$(S_i) \geq 4(d+1)$, $1 \leq i \leq k$.

In [Yu97] the bound on $k$ is not given in the nonorientable case, but from the proof one can easily derive $g/2 \leq k \leq 2g - 1$. A modification as in [T93] and the proof below gives $g/2 \leq k \leq g$. The related theorem in [Y97, Thm. 4.3] claims to prove more; instead we derive (iii) from Prop. 1, correcting a small flaw in the proofs of [T93, Y97].

Proof. The goal is to find a set of planarizing cycles $\{C_1, ..., C_k\}$ with $k = g/2$ when $\Sigma$ is orientable and with $g/2 \leq k \leq g$ when $\Sigma$ is nonorientable, with the properties stated in this theorem. We have for $\Sigma$ orientable and $g \geq 2$, $w(G, \Pi) \geq 8(d+1)(2^{g/2} - 1) \geq 8(d+1)$, and for $\Sigma'$ nonorientable and $g \geq 1$, $w(G, \Pi) \geq 8(d+1)(2^g - 1) \geq 8(d+1)$.

In all cases, let $C$ be a shortest noncontractible cycle in the embedded graph $G$; $C$ is necessarily chordless by Lemma 1. By Prop. 1, ncdist$(C) \geq 4(d+1)$ in all cases.

We cut along $C$, as chosen above, to create $G'$ with triangulation $\Pi'$ on $\Sigma'$. By Lemma 1, we have

for $g \geq 2$, $w(G', \Pi') \geq 4(d+1)(2^{g/2} - 1) \geq 4(d+1)$ when $\Sigma$ is orientable, and (*)

for $g \geq 1$, $w(G', \Pi) \geq 4(d+1)(2^g - 1) \geq 4(d+1)$ otherwise. (*)

The proof is by induction on the Euler genus $g \geq 1$, and the inductive hypothesis will be applied to the graph $G_0$ with embedding $\Pi_0$ on $\Sigma'$, defined below.

First suppose that $C$ is 2-sided and $\Sigma'$ is the sphere. Then $g = 2$, $\Sigma$ is either the torus or the Klein bottle, and we set $C_1 = C$. We have in both cases a planarizing set of one cycle $C_1 = C$ ($k = g/2$) for which ncdist$(C_1) \geq 4(d+1)$.

Next suppose that $C$ is 1-sided and $\Sigma'$ is the sphere so that $g = 1$, $\Sigma$ is the projective plane, and we set $C_1 = C$. Again we have a planarizing set of one cycle $C_1$ ($k = g$) for which ncdist$(C_1) \geq 4(d+1)$.

Thus we assume that after cutting along $C$ in $G$ on $\Sigma$, we obtain $G'$ on $\Sigma'$ which is not the sphere so that either $\Sigma$ is the Klein bottle with $g = 2$ and cutting along $C$ leaves either one or two projective planes, or else $g \geq 3$. With the bounds in (*), we have $w(G', \Pi') \geq 12(d+1)$ in all these cases. Then we consider the canonical cycles at $x'$ in $G'$ (and at $x$ if $C$ is 2-sided). Let $D'$ be the $(d+2)$-canonical cycle of $x'$ in $G'$ (and if $C$ is 2-sided, $D$ the $(d+2)$-canonical cycle of $x$ in $G'$). Also let $E'$ be the $(d+1)$-canonical cycle of $x'$ (and if $C$ is 2-sided, $E$ the $(d+1)$-canonical cycle of $x$ in $G'$). Since $w(G', \Pi') \geq 4(d+1)$, when $C$ is 2-sided, the cycles $D'$ and $D$ are disjoint and are joined by no edge. When $C$ is 1-sided, by Prop. 2 $D'$ is chordless. We delete the interior of $D'$ that contains $x'$ in $G'$ (and the interior of $D$ that contains $x$ in $G'$) and add a vertex $y'$ adjacent to all vertices of $D'$ (and add vertex $y$ adjacent to all vertices of $D$) to form $G_0$ with triangulation $\Pi_0$ on $\Sigma'$.

Now the proof will proceed by induction on the Euler genus $g \geq 2$. Since as shown in [T93, Y97] and using the inequalities in (*), $w(G_0, \Pi_0) \geq w(G', \Sigma') - 4(d+1) \geq 4(d+1)(2^{g/2} - 2) = 8(d+1)(2^{g/2-1} - 1)$ when $\Sigma'$ is orientable with $g \geq 4$, and similarly $w(G_0, \Pi_0) \geq w(G', \Sigma') - 4(d+1) \geq 8(d+1)(2^{g-1} - 1)$ when $\Sigma'$ is nonorientable with $g \geq 2$, the inductive hypothesis holds. Those inequalities also give $w(G_0, \Sigma') \geq 8(d+1)$ in these cases.

When $C$ is nonseparating, we set $C_1 = C$, we look for a planarizing set for $G_0$ on $\Sigma'$, and we add $C_1$ to that set; when $C$ is separating, we don't use $C$ in the final planarizing set. Specifically, if $C$ is 2-sided and nonseparating, $\Sigma'$ has Euler genus $g - 2$, and by induction we get a planarizing set $\{C_2, ..., C_k\}$ for $G_0$ on $\Sigma'$ with $(g-2)/2 \leq k - 1 \leq g - 2$ so that $\{C_1, ..., C_k\}$ satisfies $g/2 \leq k \leq g - 1 < g$. And when $\Sigma$ and thus $\Sigma'$ are orientable, we get $k = g/2$. Similarly when $C$ is 1-sided (so necessarily nonseparating and $\Sigma$ is nonorientable), $\Sigma'$ has Euler genus $g - 1$. By induction we get a planarizing set $\{C_2, ..., C_k\}$ for $G_0$ on $\Sigma'$ with $(g-1)/2 \leq k - 1 \leq g - 1$ so that $\{C_1, ..., C_k\}$ satisfies $g/2 < (g+1)/2 \leq k \leq g$. In both cases when $C$ is nonseparating we have that the cycles $\{C_1, ..., C_k\}$ are all chordless. By induction dist$_{G_0}(C_i, C_j) \geq d$, $2 \leq i < j \leq k$, and ncdist$_{G_0}(C_i) \geq 4(d+1)$, $2 \leq i \leq k$.

When $C$ is separating and necessarily 2-sided, we have $G_0 = G_1 \cup G_2$, $G_i$ on surface $\Sigma_i$ of Euler genus $g_i \leq g - 1$, $i = 1, 2$, with $g_1 + g_2 = g$; also $w(G_i, \Sigma') \geq w(G_0, \Sigma') \geq 8(d+1)$, $i = 1, 2$. By induction we get a planarizing set $\{C_1, ..., C_m\}$ for $G_1$, $g_1/2 \le m \le g_1$ with $m = g_1/2$ when $\Sigma_1$ is orientable, and a planarizing set $\{C_{m+1}, ..., C_k\}$ for $G_2$, $g_2/2 \le k - m \le g_2$ with $k - m = g_2/2$ when $\Sigma_2$ is orientable. Thus the planarizing set for G $\{C_1, ..., C_k\}$ has size bounded by $g/2 \le k \le g$, and since $\Sigma$ is orientable if and only if $\Sigma_1$ and $\Sigma_2$ are, then $k = g/2$. These cycles are all chordless and satisfy $\text{dist}_{G_1}(C_j, C_{j'}) \ge d$, $1 \le j < j' \le m$, $\text{dist}_{G_2}(C_j, C_{j'}) \ge d$, $m + 1 \le j < j' \le k$, and by Prop. 1 $\text{ncdist}_{G_i}(C_j) \ge w(G_0, \Sigma')/2 \ge 4(d+1)$, $1 \le j \le k$, for $i = 1, 2$. In all cases we must still check that the distances in $G$ satisfy (ii) and (iii).

We show that the desired planarizing set is $\{C_1, ..., C_k\}$, $g/2 \le k \le g$, when no $C_i$ contains $y$ or $y'$. On the other hand when, say, $C_{i'}$ contains $y'$, we modify $C_{i'}$, as follows [T93, Y97]. Namely, if $C_{i'}$ contains the subpath $uy'v$ in $G_0$ with $u, v \in D'$, then we replace this subpath with a *shortest* path in $G$ of the form $u, y_1, ..., y_r, v$ with $y_1, ..., y_r$ a subpath of $E'$, all vertices at distance one from $D'$. We do the same when $C_{j'}$ contains $y$ (possibly $j' = i'$). With these choices the cycles $C_{i'}$ and $C_{j'}$ remain chordless (since $E'$, $E$, $C_{i'}$, and $C_{j'}$ are chordless and the subpath is chosen to be shortest). The distance between cycles $C_j$ and $C_{j^*}$, $1 \le j < j^* \le k$, given by induction on $G_i$, $i = 0, 1, 2$, remains at least $d$ on $G$ except possibly when a shortest path joining them passes through $y$ or $y'$. It is a routine check to see that the length of such a shortest path cannot decrease in $G$. Due to the construction of $G_0$ (when $C$ is nonseparating), $\text{dist}_G(C_1, C_i) \ge d$, $2 \le i \le k$. Due to the construction of $G_1$ and $G_2$ (when $C$ is separating), the distance between a cycle on $G_1$ and one on $G_2$ is at least $2d$. Thus (ii) holds. As noted above when $C$ is nonseparating, $\text{ncdist}_G(C_1) \ge 4(d+1)$. Also for $1 \le j \le k$ if $P$ is a path giving $\text{ncdist}_{G_i}(C_j)$ by induction for $i = 0, 1, 2$, then the minimum length given by $l(P)$ is not affected unless $P$ passes through $y'$ or $y$. Again in $G$ a longer alternate path must be found or else the ncdist-defining path itself does not decrease in length when modified or replaced to become an ncdist-defining path in $G$. Thus (iii) also holds. QED

**2. The nice cycle lemma**

The driving force behind our results and many of their predecessors is the *nice cycle lemma*. This says that if $C$ is a noncontractible cycle in an "orderly" triangulation (defined below), then there is a chordless cycle $C'$ that is near to $C$, homotopic to $C$, and particularly nice from the point of view of 5-coloring. A chordless noncontractible cycle $C$ in an embedded graph is said to be *nice* if it either has even length or contains a vertex of degree at most 4. The first published version of a nice cycle lemma appeared in [AS82]. There the cycle $C$ had to be a minimum length noncontractible cycle in an orderly triangulation of the torus, and the cycle $C'$ was found in the first neighborhood of $C$. Stromquist established a nice cycle lemma for an arbitrary surface where $C'$ was found in the "filled first neighborhood" of $C$ [S77], but this was never published. The filled first neighborhood might contain vertices arbitrarily far from $C$. Thomassen constructed a nice cycle in [T93]; however, to obtain his version of $C'$ he used the notion of weak geodesic, requiring a considerable detour away from $C$. Using local modifications he found his nice cycle on one side of $C$, but, as we now know, finding a nice cycle requires zigzagging to both sides of $C$. We sketch the argument at the end of this section. Finally we need a nice cycle lemma for both nonorientable and orientable surfaces; a subtle flaw was found by the third author in the proof of [AH04].

As in [AH04] we introduce "optimal *k*-shortcuts" and use this idea to give a proof of the lemma. Our result also has the advantage that all vertices of $C'$ are within distance four of the original $C$.

If $G$, $\Pi$ is a triangulation of the surface $\Sigma$ and $U \subset V(G)$, $G[U]$ is said to be *orderly* if every contractible 3-cycle in $G[U]$ bounds a face of $G$ and if every contractible 4-cycle in $G[U]$ is either the boundary of two triangles of $G$ that share an edge or the first neighbor circuit of a vertex of degree 4 in $G$. Orderly graphs have been helpful in inductive proofs of 5-coloring theorems for embedded graphs [AS82, H84].

Suppose $C = u_1, u_2, ..., u_t$, $t > 6$, is a chordless, noncontractible cycle in a graph embedded on a surface $\Sigma$ of Euler genus $g > 0$. A path of length 6, say $v_1, v_2, ..., v_7$, is said to be an *optimal 6-shortcut* for *C* if $v_1 = u_i$ for some *i*, $v_7 = u_j$ for some *j*, the resulting noncontractible cycle
$C' = u_1, \ldots, \{u_i = v_1\}, v_2, \ldots, v_6, \{v_7 = u_j\}, u_{j+1}, \ldots, u_t$ is homotopic to *C*, *and C' is as short as possible.*
(Note that by renumbering, if the shortcut goes the "other" way and includes $u_1$ in its bypass, by relabeling we can assume the shortcut is as above.) Equivalent to the first condition, the cycle
$u_i, u_{i+1}, \ldots, u_j, v_6, v_5, \ldots, v_2, u_i$ is contractible. By Thomassen's 3-path-property [T93, Prop. 4.3.1], one subpath of *C* between $u_i$ and $u_j$ together with the 6-shortcut is a noncontractible cycle; without loss of generality that cycle is *C'*. Note that $v_2 \ldots, v_6$ may be vertices of *C*, that possibly *C' = C*, and that when *C' ≠ C*, *C'* might not be chordless.

Let *C* be a noncontractible cycle in an embedded graph and *K* a chord of *C*. We say *K* is a *contractible chord* if *K* together with one subpath of *C* between the endpoints of *K* is contractible; otherwise *K* is called a *noncontractible chord*. We prove a variation on the nice cycle lemma of [AH04], correcting an error in the proof there.

Lemma 3 (Nice Cycle Lemma; see also [AH04]). For $d \geq 2$, let a simple graph *G* be a 2-cell triangulation Π of a surface Σ of positive Euler genus $g \geq 1$ with $w(G, \Pi) \geq 8(d+1)(2^g - 1)$. Let *C* be a chordless noncontractible cycle in a planarizing collection of cycles, as given in Thm. 1. Suppose $G_{C,4} = G[C \cup N(C) \cup \ldots \cup N^4(C)]$ is orderly. Then there is a chordless noncontractible cycle *C'* in $G_{C,4}$ such that *C'* is homotopic to *C* and *C'* is nice. Also $\text{ncdist}(C') \geq 4d - 2$, and for $X \subset V(G)$ disjoint from $G_{C,4}$ $\text{dist}(X, C') \geq \text{dist}(X, C) - 4$.

Proof. Suppose $C = u_1, u_2, \ldots, u_t$ is a chordless noncontractible cycle in a planarizing collection for *G* with $\text{ncdist}(C) \geq 4(d+1) \geq 12$ as given by Thm. 1; necessarily $t \geq w(G, \Pi) \geq 8(d+1)(2^g - 1) \geq 24$. Choose $M = v_1, v_2, \ldots, v_7$ to be an optimal 6-shortcut for *C* where $v_1 = u_i$ and $v_7 = u_j$. From the definition of optimal shortcut, the labels may be chosen so that $1 \leq i < j \leq t$ and $M \cup u_i, u_{i+1}, \ldots, u_j$ is contractible. Setting $C' = u_j, u_{j+1}, \ldots, u_t, u_1, \ldots, u_i \cup M$, we have that *C'* is homotopic to *C*; in particular, it is noncontractible. It might be that *C' = C* (for example when *C* is a shortest noncontractible cycle in *G*, Π), and then *C'* is chordless and $\text{ncdist}(C') \geq 4(d+1)$.

We claim that *C'* has no chord also when *C' ≠ C*. Suppose *K* were a chord of *C'*. Since *C* is chordless, at least one endpoint of *K* is not in *C* and therefore is in *M*. If *K* is contractible, then *M* is not an optimal shortcut, as *C'* can be shortened; therefore we may assume *K* is noncontractible.

Using *M* to join the endpoints of *K* to *C* by disjoint paths, we get a path *P* of length at most 7 joining two vertices of *C*. Since *K* is noncontractible for *C'*, none of the cycles in $C \cup P$ is contractible, showing that $\text{ncdist}(C) \leq 7$, contradicting the fact that $\text{ncdist}(C) \geq 12$.

We also claim that $\text{ncdist}(C') \geq 4d$ when *C' ≠ C*. Let *P* be a path defining $\text{ncdist}(C')$; if it is a cycle then $\text{ncdist}(C') \geq w(G, \Pi) \geq 8(d+1)$. If *P* joins two vertices of *C'\M*, then $l(P) \geq 4(d+1)$. Otherwise *P* joins two vertices of *M* or a vertex of *M* with one in *C'\M*, and *P* together with each subpath of *C'* is noncontractible. Then *P* plus at most 4 edges is either a noncontractible cycle or forms two noncontractible cycles with subpaths of *C*, and we conclude that $l(P) + 4 \geq \text{ncdist}(C) \geq 4(d+1)$ so that $l(P) \geq 4d$.

Since $w(G, \Pi) \geq 8(d+1)(2^g - 1) \geq 24$, $R^4(M) \cap L^4(M) = \emptyset$ by Cor. 1 and Prop. 2, and no edge joins $R^4(M)$ with $L^4(M)$. We may assume that both *C* and *C'* have odd length and contain no vertex of degree at most 4. Thus both $t (= |C|)$ and $i - 1 + 7 + t - j (= |C'|)$ are odd. If *D* denotes the subpath of *C'* from $u_{j+1}$ to $u_{i-1}$ inclusive, then $|D| = i - 1 + t - j$ contains an even number of vertices.

*Case (i)*: Suppose $\deg_R(v_3) \geq 2$ and $\deg_R(v_4) \geq 2$. Let *x* denote the vertex that is in a triangle to the right of the edge joining $v_3$ with $v_4$. Suppose $v_1, A, x$ is a shortest path from $v_1$ to *x* among the vertices in $N(v_3)$ (with *A* not including either endpoint). Normally *A* includes $v_2$ but that is not necessary. Next let $x, B, v_6$ be a shortest path from *x* to $v_6$ among the neighbors of $v_4$ (with *B* not including either endpoint.)

Normally $B$ includes $v_5$ but that is not necessary. Note that $v_4$ or $x$ might be at distance 3 from $C$ and perhaps a neighbor of $v_4$ is at distance 4 from $C$. Otherwise all vertices are at distance 3 or less from $C$.

Consider the cycle $C^{(2)} = D, v_1, A, x, v_4, v_5, v_6, v_7$. $C^{(2)}$ is homotopic to $C'$ since ncdist($C'$) $\geq 4d \geq 8$ and all vertices of $C^{(2)}$ are at distance at most 1 from $C'$. Thus ncdist$(C^{(2)}) \geq 4d - 2$. Let's examine the possible chords for $C^{(2)}$. From the latter inequality we see that a noncontractible chord, with one end in $A$ or at $x$ plus at most two edges, joins two vertices of $C'$, a contradiction. And if there were a contractible chord from $D$ to a vertex in the subpath $A, x, v_4, v_5, v_6$, then we didn't select our optimal 6-shortcut correctly. This is also the case if there were a contractible chord from $A$ to either $v_6$ or $v_7$. If there were a contractible chord from $A$ to $v_4$, then $G_{C,4}$ contains a contractible 3-cycle that is not a face boundary. If there were a contractible chord from $A$ to $v_5$, then we have a contractible 4-cycle, say $a, v_3, v_4, v_5$, with $x$ in its interior for some $a \in A$. If $x$ is adjacent to $v_5$ and deg$(x) = 4$, then $D, v_1, v_2, v_3, x, v_5, v_6, v_7$ is nice. If $x$ is not adjacent to $v_5$ or deg$(x) \neq 4$, then $G_{C,4}$ contains a forbidden contractible 4-cycle. Thus either the number of vertices in $A$, $|A|$, is odd or we are done since $C^{(2)}$ would be nice.

Next consider the cycle $C^{(3)} = D, v_1, v_2, v_3, x, B, v_6, v_7$ with all vertices at distance at most 4 from $C$ and at distance at most 1 from $C'$. As in the preceding paragraph either $|B|$ is odd or we are done since $C^{(3)}$ would be nice. Finally, look at $C^{(4)} = D, v_1, A, x, B, v_6, v_7$. $C^{(4)}$ is homotopic to $C'$ and $|C^{(4)}|$ is even. The only chord we have not considered is an edge from $A$ to $B$. If there were a noncontractible chord joining $A$ and $B$, then the chord plus two edges gives ncdist($C'$) $\leq 3$, a contradiction. If there were a contractible chord joining $A$ and $B$, then it produces a forbidden contractible 4-cycle. In addition we have, for $i = 3, 4$, ncdist$(C^{(i)}) \geq 4d - 2$ since if $P$ is a defining path for ncdist$(C^{(i)})$, then either $l(P) \geq 4d$ or $l(P) + 2 \geq 4d$.

*Case (ii)*. Suppose on $C'$ we have deg$_R(v_3) < 2$. If either deg$_R(v_3) = 0$ or deg$_R(v_3) = $ deg$_R(v_4) = 1$, then $v_1, ..., v_7$ is not an optimal 6-shortcut. Thus we assume that deg$_R(v_3) = 1$ and deg$_R(v_4) \geq 2$. Symmetrically we may assume deg$_R(v_5) = 1$. If either deg$_L(v_5) = 1$ or deg$_L(v_3) = 1$, then $C'$ is nice. Now if deg$_L(v_4) > 1$, then we are back in Case (i) interchanging left with right. Thus we assume deg$_R(v_3) = $ deg$_L(v_4) = 1$ and deg$_L(v_5) \geq 2$. Suppose $x$ is $v_3$'s unique right neighbor. We alter $C'$ by replacing $v_3$ with $x$. In the resulting cycle both deg$_L(v_4)$ and deg$_L(v_5)$ are at least 2 and we are in Case (i) again.

*Case (iii)*. Suppose on $C'$ we have deg$_R(v_3) \geq 2$, and deg$_R(v_4) < 2$. As in Case (ii) we may assume deg$_R(v_4) = 1$. If deg$_R(v_5) \leq 1$, then $v_1, ..., v_7$ is not an optimal 6-shortcut so that we have deg$_R(v_5) \geq 2$. Now deg$_L(v_4) \geq 2$ since otherwise $C'$ is not an optimal 6-shortcut or else it has a vertex of degree 4 and is nice. Similarly and by Case (i) we must have deg$_L(v_3) = $ deg$_L(v_5) = 1$. Then by reversing the direction of the 6-shortcut we are in Case (ii) again.

Thus in each case we find the nice cycle $C'$ within $G_{C,4}$ with ncdist($C'$) $\geq 4d - 2$. Each nice cycle is no more than four edges away from its original planarizing cycle. Thus for $X$ outside of $G_{C,4}$ dist$(X, C') \geq$ dist$(X, C) - 4$. QED

A closer reading of the proof shows that the results hold with $d \geq 2$ when $g = 1$ and $d \geq 1$ otherwise. Also we note that when $\Sigma$ is orientable, the Nice Cycle Lemma holds when $w(G, \Pi) \geq 8(d+1)(2^{g/2} - 1)$ and $d \geq 2$, but for simplicity we assume uniform bounds.

We digress a moment to describe how to construct a noncontractible cycle $C$ in a graph embedded in a nonorientable surface with the property that a nearby nice cycle must use both sides of $C$. Let $C$ consist of vertices whose right degrees are alternately 1 and 3. It is easy to build such a cycle in a 6-regular triangulation of a Klein bottle, and from there it can be on any nonorientable surface. For such a cycle, there is no nice replacement cycle lying locally on one side of $C$. Thus when $C$ has odd length and is 1-sided, the nice cycle must be found using detours locally to both sides of *C*.

### 3. A 5-color theorem for locally planar graphs

A version of the next theorem for orientable surfaces appeared in [AH01]. There are several good reasons to have another look. First, the proof given below works for both orientable and nonorientable surfaces. Second, the particular nice cycle lemma needed in the earlier proof never appeared in print. Third the constants are better.

**Theorem 2 (**See also [AH04]) Suppose *G* is an embedded triangulation Π on a surface Σ of Euler genus $g > 0$, $d \geq 14$, and $w(G, \Pi) \geq 8(d+1)(2^g - 1)$. Then *G* is 5-colorable.

Proof. We apply Theorem 1 with $d = 14$ to get a planarizing set of cycles $C_1, C_2, ..., C_s$, $g/2 \leq s \leq g$, with lower bounds on the distance between pairs of these cycles and on ncdist($C_i$), $1 \leq i \leq s$.

First, for each $i \leq s$ we make the subgraph $G_{C_i,4}$ (see Lemma 3). For $i \leq s$ these subgraphs are disjoint and joined by no edge by Cor. 1, Prop. 2, and Thm. 1 since $d = 14$. For each $G_{C,4}$ we delete vertices of *G* interior to every contractible 3-cycle of $G_{C,4}$ and, if there is more than one, all vertices interior to every contractible 4-cycle of $G_{C,4}$. In the latter case we add a vertex adjacent to all four boundary vertices. Call this intermediate graph $G_1$. Then for $1 \leq i \leq s$ we apply Lemma 3 with $C = C_i$ and let $C'_i$ denote the resulting nice cycle, homotopic to $C_i$, lying within $G_{C_i,4}$. This gives a new planarizing set $\{C'_1, ..., C'_s\}$.

By Lemma 3 we have dist($C'_i, C'_j$) $\geq d - 8 \geq 6$ for $i \neq j$ and ncdist($C'_i$) $\geq 4d - 2 \geq 54$ for $1 \leq i \leq s$. With this planarizing set for $G_1$ we form a planar triangulation as follows (in a manner slightly different than "cutting" along these cycles). For each $C'_i$ that is 2-sided, $1 \leq i \leq s$, we remove $C'_i$ and replace it with two vertices $x_i$ and $y_i$. Then $x_i$ is made adjacent to each vertex of $L(C'_i)$ and $y_i$ is adjacent to each vertex of $R(C'_i)$. For $C'_i$ that is 1-sided, we remove $C'_i$ and replace it with $x_i$ adjacent all vertices in $N(C'_i)$. (It does not matter that some multiple edges might be introduced.) The resulting graph $G_0$ is a plane triangulation with the distance between every pair of vertices from $\{x_i, y_i, x_j, y_j\}$ at least $d - 8 \geq 6$ for $1 \leq i \leq j \leq s$ (and the same distance bounds hold when some $y_i$ or $y_j$ does not exist due to a 1-sided cycle.)

Suppose *c* is a 4-coloring of $G_0$ using colors {1,2,3,4} [ApHa76, RSST97]. For each $i \leq s$ such that $c(x_i) \neq c(y_i)$, we recolor each vertex in $N^2(y_i)$ that is colored $c(y_i)$ with color 5 and then at $y_i$ we perform a $(c(x_i), c(y_i))$-Kempe change so that $y_i$ gets the same color as $x_i$. Since the distance between the pairs of $\{x_i, x_j, y_i, y_j\}$, $i \leq j$, is at least 6, this recoloring is valid. Next we transfer this coloring back to $G_1$. For each *i* for which $C'_i$ has even length, we 2-color $C'_i$ with $\{c(x_i), 5\}$, and for each $C'_i$ with a vertex $v_i$ of degree 4, we 2-color $C'_i$ with $\{c(x_i), 5\}$ except for $v_i$ which is colored last with whatever color is available. Call this coloring of $G_1$ $c_1$.

Then we add back and color the vertices deleted within contractible 3- and 4-cycles lying within $G_{C_i,4}$ for each $i \leq s$. Let *C* be a contractible 3-cycle in $G_{C_i,4}$ with nonempty interior in *G*. Let *H* be the plane graph Int(*C*). We can easily 4-color *H* so that it agrees with the coloring of *C* from $G_1$.

Suppose $C = u, v, x, y$ is a contractible 4-cycle within $G_{C_i,4}$ that within *G* contains more than one interior vertex so that some interior vertices were deleted to form $G_1$. Let *H* be the plane graph Int(*C*). We may assume that in the coloring of $G_1$, $c_1(u) = 1$ and $c_1(v) = 2$.

*Case* (i). Suppose $c_1(x) = 1$. We transform the plane graph *H* by identifying the vertices *u* and *x* within the outer face of *H*. If $c_1(y) = 2$ (resp., 3), then we 4-color *H* using colors {1,3,4,5} (resp., {1,2,4,5}) making sure that *u* and *x* (resp., *u*, *v*, and *x*) are correctly colored. Then we assign color 2 to both *v* and *y* (resp., color 3 to *y*) for a coloring of *H* that agrees on the boundary with that of *C* in $G_1$.

*Case* (ii). Suppose $c_1(x) = 3$. We transform the plane graph *H* by adding an edge joining *u* and *x* in the outerface of *H*. If $c_1(y) = 2$ (resp., 4), we 4-color *H* with {1,3,4,5} (resp., {1,2,3,5}) making sure *u* and *x* (resp., *u*, *v*, and *x*) are correctly colored. We assign color 2 to *v* and *y* (resp., color 4 to *y*). This gives a 5-coloring of *G*. QED

We note again that when Σ is orientable, the bound $w(G, \Pi) \geq 8(d+1)(2^{g/2} - 1)$ is sufficient with $d \geq 14$.